# Roots of any Polynomial with Complex Integer Coefficients using Replacement Sequences, Ruler and Compass

Ashok Kumar Mittal,
Department of Physics,
Allahabad University, Allahabad – 211 002, India
(**Email address:** mittal_a@vsnl.com)

Ashok Kumar Gupta,
Department of Electronics and Communication,
Allahabad University, Allahabad - 211 002, India
(Email address: akgjkiapt@hotmail.com)

**Abstract:** The roots of any polynomial of degree *m* with complex integer coefficients can be computed by manipulation of sequences made from distinct symbols and counting the different symbols in the sequences. This method requires only 'primitive' operations like replacement of sequences and counting of symbols. No calculations using 'advanced' operations like multiplication, division, logarithms etc. are needed. The method can be implemented as a geometric construction using only a ruler and a compass.

## 1. Introduction

Rational approximates to the largest (absolute value) real root of any polynomial of degree *m* with integer coefficients can be obtained to arbitrary accuracy by a method requiring only operations consisting of replacement of sequences and counting of symbols [5]. The replacement rules for this method are obtained from a matrix, called the 'replacement matrix', whose elements are related to the coefficients of the polynomial. The method can be extended to obtain complex roots by modifying the 'replacement matrix'.

The eigenvalues of the 'replacement matrix' are the roots of the polynomial multiplied by the coefficient of the highest power of the polynomial. The largest eigenvalue can be obtained by the method of matrix iteration, provided the largest eigenvalue is unique [1]. It is possible to carry out this matrix iteration method using only 'primitive operations' like replacement of sequences and counting of symbols even to obtain complex roots by appropriately modifying the 'replacement matrix'.

The method can be implemented geometrically using only a ruler and a compass.

## 2. Replacement Matrix and the Replacement Rules

Let $p(x) = a_0 x^m + a_1 x^{m-1} + \ldots + a_m$, where $a_i$ are integers, be the given polynomial with integer coefficients. The 'replacement matrix' corresponding to the polynomial is defined by

$$\mathbf{R} = \begin{bmatrix} -a_1 & -a_2 & -a_3 & \ldots & -a_{m-1} & -a_m \\ a_0 & 0 & 0 & \ldots & 0 & 0 \\ 0 & a_0 & 0 & \ldots & 0 & 0 \\ . & . & . & \ldots & . & . \\ . & . & . & \ldots & 0 & . \\ 0 & 0 & 0 & \ldots & a_0 & 0 \end{bmatrix} \qquad (1)$$

Let $r_i$, $i = 1, 2, \ldots, m$ be the roots of $p(x)$. Then the eigenvalues of $\mathbf{R}$ are $a_0 r_i$ with corresponding eigenvector $[r_i^{m-1}, r_i^{m-2}, \ldots, r_i, 1]^T$.

Let $A = \{A_1, A_2, \ldots, A_m\}$ be a finite set of $m$ symbols. Elements of A may be called 'letters' belonging to the 'alphabet' A. Let $A^p = A \times A \times \ldots \times A$ (p times). Elements of $A^p$ are called 'words' of length p that can be formed from 'letters' of A. Let $A^* = \bigcup_p A^p$. Clearly, $A^*$ consists of all 'words' that can be made from the 'alphabet' A. The replacement rule, that is used to manipulate the symbolic sequences consisting of these symbols, is obtained from the 'replacement matrix' $\mathbf{R}$.

Let us first consider the case when all the elements of this matrix are non-negative. Let $R_{ij}$ denotes the element in the $i^{th}$ row and $j^{th}$ column of the matrix $\mathbf{R}$. The replacement rule R, gives a rule to replace a word W by a new word $W^*$. This rule consists of replacing each letter $A_j$ in W by a word consisting of the $i^{th}$ symbol repeated $R_{ij}$ times. Thus,

$$R(A_j) = (A_1)^{\wedge}R_{1j} \ (A_2)^{\wedge}R_{2j} \ldots \ldots (A_i)^{\wedge}R_{ij} \ldots \ldots (A_m)^{\wedge}R_{mj} \qquad (2)$$

where $(A_i)^{\wedge}k = A_i^k$ stands for repeating k times the symbol $A_i$. This replacement rule fails to have any meaning if any element of the matrix $\mathbf{R}$ is negative.

In order to find a replacement rule applicable for the case when some of the elements of the 'replacement matrix' are negative integers, we enlarge the 'alphabet' A to have $2m$ symbols. We denote this enlarged 'alphabet by $A = \{A_1, A_2, \ldots, A_m, A_{m+1}, \ldots, A_{2m}\}$. The symbols $A_{m+i}$ are also denoted by $\tilde{A}_i$. The symbols $A_i$ and $\tilde{A}_i$ will be called conjugate to each other, so that $A_i^c = \tilde{A}_i$ and $(\tilde{A}_i)^c = A_i$. Then $(A_i)^{-k}$ stands for repeating k times the symbol $A^c$. The replacement rule (2) will have to be extended to provide for replacement of $2m$ symbols. This is given by

$$R(A_j) = (A_1)^{\wedge}R'_{1j} \ (A_2)^{\wedge}R'_{2j} \ldots \ldots (A_i)^{\wedge}R'_{ij} \ldots \ldots (A_{2m})^{\wedge}R'_{(2m)j} \qquad (3)$$

where $R_{ij}$ are now matrix elements of the $2m \times 2m$ extended matrix

$$\mathbf{R}_e = \begin{bmatrix} \mathbf{R} & 0 \\ 0 & \mathbf{R} \end{bmatrix} \qquad (4)$$

Let $W = A_{s1}A_{s2}\ldots A_{sq}$, where $A_{sk}$ $(k=1,2,\ldots,q) \in A$. Then $W \in A^*$. The replacement rule R induces a mapping $R^*: A^* \to A^*$ defined by

$$W^* = R^*(W)$$
$$= R^*(A_{s1}A_{s2}\ldots\ldots A_{sq})$$
$$= R(A_{s1})R(A_{s2})\ldots\ldots R(A_{sq}) \qquad (5)$$

Let $W_0 \in A^*$. Then

$$W_i = (R^*)^i (W_0) = R^*((R^*)^{i-1}(W_0)) = R^*(W_{i-1}) \qquad (6)$$

denotes the word obtained by i times repeated application of the replacement rule (3) on the initial word $W_0$.

Let $\mathbf{n}: A^* \to V^m$ be a mapping which assigns to a word in $A^*$ a vector

$$\mathbf{n}(W) = (n_1(W) - \tilde{n}_1(W), \ldots, n_i(W) - \tilde{n}_i(W), \ldots, n_m(W) - \tilde{n}_m(W))^T \qquad (7)$$

in an m-dimensional vector space $V^m$ such that $n_i(W)$ and $\tilde{n}_i(W)$ are non-negative integers denoting the number of times that $A^+_j$ and $\tilde{A}_i$ occur, respectively, in the word W. It is straightforward to see that

$$\mathbf{n}(W^*) = \mathbf{n}(R^*(W)) = \mathbf{R}(\mathbf{n}(W)) \qquad (8)$$

Equations (6) and (8) imply

$$\mathbf{n}(W_i) = \mathbf{n}(R^*(W_{i-1})) = \mathbf{R}(\mathbf{n}(W_{i-1})) = \mathbf{R}^i(\mathbf{n}(W_0)) \qquad (9)$$

For almost any $W_0$, $\mathbf{n}(W_i)$ tends, as $i \to \infty$, to multiples of the eigenvector of $\mathbf{R}$ corresponding to the eigenvalue $\lambda$ having maximum absolute value [1], provided it is unique. Thus one finds that for almost any initial word $W_0$,

$$\lim_{i \to \infty} n_j(W_i)/n_{j+1}(W_i) = \lambda \qquad j = 1, 2, \ldots, m-1 \qquad (10)$$

The operation of subtraction in equation (7) can be carried out by removing pairs of symbols $A_i$ and $\tilde{A}_i$ from each replacement sequence. Thus, Eqn (10) shows that the largest real root of any polynomial with integer coefficients can be obtained to any accuracy only by operations consisting of replacement of sequences and counting of symbols.

The method is made more versatile by noting that the eigenvalues of the matrix $\mathbf{R}' = (\alpha \mathbf{I} + \beta \mathbf{R})$ are $r_i' = (\alpha + \beta a_0 r_i)$ with corresponding eigenvectors $[r_i^{m-1}, r_i^{m-2}, \ldots,$

$r_i$, $1]^T$, where **I** is the *m x m* identity matrix. A different set of replacement rules may be obtained by replacing matrix **R** in (1) by **R'**

The method would seem to be incapable of obtaining complex roots. However, to obtain complex roots it is possible to extend the method by letting $\alpha$ and $\beta$ be complex integers and replacing each complex number $a + ib$ in the replacement matrix **R'** by a 2 x 2 matrix

$$\begin{vmatrix} a & b \\ -b & a \end{vmatrix}$$

Let $\mathbf{R}_c$ denote the matrix so obtained. Let

$$\mathbf{w}^0 = [u^0_1, v^0_1, u^0_2, v^0_2, \ldots, u^0_m, v^0_m]^T \qquad (11)$$

be an arbitrary real *2m* dimensional vector with integer elements. Let

$$\mathbf{R'}^k \mathbf{w}^0 = [u^k_1, v^k_1, u^k_2, v^k_2, \ldots, u^k_m, v^k_m]^T \qquad (12)$$

Then,

$$\lim_{k \to \infty} (u^k_j - iv^k_j)/(u^k_{j+1} - iv^k_{j+1}) = r_{max}, \quad j = 1, 2, \ldots, m-1 \qquad (13)$$

where $r_{max}$ is the root with the largest absolute value provided it is unique. The integers $u^k_j$, $v^k_j$ can be obtained using the replacement rules corresponding to the matrix $\mathbf{R}_c$, appropriately extended, if there are negative integers in it. It is only at the last stage, in equation (13), that just one division of two complex numbers is required. This operation can be carried out geometrically, using a straight edge and a compass only as described in the next section.

In this way even complex roots can be obtained to arbitrary accuracy without requiring any 'advanced' operation such as multiplication, division, logarithms and anti-logarithms, although the justification of the method uses 'advanced' concepts like the theory of matrices.

### 3. Geometrical implementation of the method

One can implement the above method as a geometrical construction of a line segment in the complex plane that represents a complex root of the polynomial. This construction uses only the straight edge, a compass, and pencils of at most *4m* different colors $c_i$ in one to one correspondence with the symbols used in the replacement rules. One begins by drawing an arbitrary line made from unit segments of different colors. Below this line, another line is drawn, in which each unit segment of the first line is replaced by a sequence of unit segments in accordance with the replacement rules. This process is continued an arbitrary number of times. Thereafter, on the next line one first draws a unit

$c_1$-colored segment for every unit $c_1$-colored segment in the last line one after another. This is followed by $c_2$, $c_3$ and $c_4$ -colored segments in the same way.

In the line so drawn, we mark the beginning of the $c_1$ line as the point A. The junction point between $c_1$ and $c_2$ is marked as B, that between $c_2$ and $c_3$ as C, between $c_3$ and $c_4$ as D. The end of the $c_4$ line is marked as E. The line AE is placed along OX, the positive real axis of the complex plane, with the point A at the origin O. The segment BC is rotated clockwise by a right angle to obtain the line BC'. AC' is joined. On the positive real axis of the complex plane a segment AD' equal to CD and D'E" are marked off. The segment D'E" is rotated clockwise by a right angle to obtain the line D'E'. AE' is joined. The line AE' is rotated to AE" where E" lies on OX. The line AC' is rotated by an angle equal to angle E'AE" to become AC". On OX a unit segment AU is marked off. E"C" is joined. From U a line parallel to E"C" is drawn to intersect AC" at R. Then AR represents a complex number in the complex plane that approximates a complex root of the given polynomial. All the above operations can be carried out using a ruler and compass.

Arbitrary accuracy can be obtained by increasing the number of times that the replacement rule is applied.

## 4. Example

Let

$$p(x) = x^2 + 1 \tag{14}$$

be the given polynomial. Here $m = 2$, $a_0 = 1$, $a_1 = 0$ and $a_2 = 1$. The 'replacement matrix' corresponding to this polynomial is given by

$$\mathbf{R} = \begin{bmatrix} 0 & -1 \\ 1 & 0 \end{bmatrix} \tag{15}$$

Choosing $\alpha = i$ and $\beta = 1$, we get,

$$\mathbf{R}' = \begin{bmatrix} i & -1 \\ 1 & i \end{bmatrix} \tag{16}$$

Replacing each number of the form $a + ib$ in (16) by a 2 x 2 matrix

$$\begin{vmatrix} a & b \\ -b & a \end{vmatrix}$$

we get,

$$\mathbf{R}_c = \begin{bmatrix} 0 & 1 & -1 & 0 \\ -1 & 0 & 0 & -1 \\ 1 & 0 & 0 & 1 \\ 0 & 1 & -1 & 0 \end{bmatrix} \tag{17}$$

In order to deal with negative integers, the extended replacement matrix is given by

$$\mathbf{R}_e = \begin{bmatrix} \mathbf{R}_c & 0 \\ 0 & \mathbf{R}_c \end{bmatrix} \tag{18}$$

Using the alphabet $A = \{0,1,2,3,\tilde{0},\tilde{1},\tilde{2},\tilde{3}\}$, this 'replacement matrix' gives rise to the following replacement rules:

$$\begin{aligned} 0 &\to \tilde{1}2 \\ 1 &\to 03 \\ 2 &\to \tilde{0}\tilde{3} \\ 3 &\to \tilde{1}2 \\ \tilde{0} &\to 1\tilde{2} \\ \tilde{1} &\to \tilde{0}\tilde{3} \\ \tilde{2} &\to 03 \\ \tilde{3} &\to 1\tilde{2} \end{aligned} \tag{19}$$

This replacement rule generates the following sequence of sequences:

0

$\tilde{1}2$

$\tilde{0}\tilde{3}\tilde{0}\tilde{3}$

$1\tilde{2}1\tilde{2}1\tilde{2}1\tilde{2}$

0303030303030303

The number count of different symbols in the above sequences can be summarized in the following table:

| **n(0)** | **n(1)** | **n(2)** | **n(3)** | **(n(0) - in(1))/(n(2) - in(3))** |
|---|---|---|---|---|
| 1 | 0 | 0 | 0 | |
| 0 | -1 | 1 | 0 | i |
| -2 | 0 | 0 | -2 | i |
| 0 | 4 | -4 | 0 | i |
| 8 | 0 | 0 | 8 | i |

As expected the ratio $(\mathbf{n}(0) - i\mathbf{n}(1))/(\mathbf{n}(2) - i\mathbf{n}(3))$ converges to the largest (absolute value) root, namely i of polynomial p(x) in (14). Choosing $\mathbf{R}' = -i\mathbf{I} + \mathbf{R}$ would yield the other root, namely -i.

## 5. Another example

Let

$$p(x) = x^2 - i \tag{20}$$

be the given polynomial. Here $m = 2$, $a_0 = 1$, $a_1 = 0$ and $a_2 = -i$. The 'replacement matrix' corresponding to this polynomial is given by

$$\mathbf{R} = \begin{bmatrix} 0 & i \\ 1 & 0 \end{bmatrix} \tag{21}$$

Choosing $\alpha = i$ and $\beta = 1$, we get,

$$\mathbf{R}' = \begin{bmatrix} i & i \\ 1 & i \end{bmatrix} \tag{22}$$

Replacing each number of the form a + ib in (22) by a 2 x 2 matrix

$$\begin{vmatrix} a & b \\ -b & a \end{vmatrix}$$

we get,

$$\mathbf{R}_c = \begin{bmatrix} 0 & 1 & 0 & 1 \\ -1 & 0 & -1 & 0 \\ 1 & 0 & 0 & 1 \\ 0 & 1 & -1 & 0 \end{bmatrix} \quad (23)$$

In order to deal with negative integers, the extended replacement matrix is given by

$$\mathbf{R}_e = \begin{bmatrix} \mathbf{R}_c & \mathbf{0} \\ \mathbf{0} & \mathbf{R}_c \end{bmatrix} \quad (24)$$

Using the alphabet A = {0,1,2,3,0˜,1˜,2˜,3˜}, this 'replacement matrix' gives rise to the following replacement rules:

$$\begin{aligned}
0 &\to 1\tilde{2} \\
1 &\to 03 \\
2 &\to 1\tilde{3}\tilde{} \\
3 &\to 02 \\
0\tilde{} &\to 1\tilde{2}\tilde{} \\
1\tilde{} &\to 0\tilde{3}\tilde{} \\
2\tilde{} &\to 13 \\
3\tilde{} &\to 0\tilde{2}\tilde{}
\end{aligned} \quad (25)$$

This replacement rule generates the following sequence of sequences:

0

1˜2

0˜3˜1˜3˜

12˜0˜2˜0˜3˜0˜2˜

031312˜1312˜0˜2˜12˜13

The number count of different symbols in the above sequences can be summarized in the following table:

| **n**(0) | **n**(1) | **n**(2) | **n**(3) | (**n**(0) - i**n**(1))/(**n**(2) - i**n**(3)) |
|---|---|---|---|---|
| 1 | 0 | 0 | 0 | inf |
| 0 | -1 | 1 | 0 | i |
| -1 | -1 | 0 | 2 | 0.5000 + 0.5000 i |
| -3 | 1 | -3 | -1 | 0.8000 + 0.6000 i |
| 0 | 6 | -4 | 4 | 0.7500 + 0.7500 i |
| 10 | 4 | 4 | 10 | 0.6897 + 0.7241 i |
| 14 | -14 | 20 | 0 | 0.7000 + 0.7000 i |
| -14 | -34 | 14 | -34 | 0.7101 + 0.7041 i |
| -68 | 0 | -48 | -48 | 0.7083 + 0.7083 i |
| -48 | 116 | -116 | 48 | 0.7066 + 0.7076 i |
| 164 | 164 | 0 | 232 | 0.7069 + 0.7069 i |
| 396 | -164 | 396 | 164 | 0.7072 + 0.7070 i |
| 0 | -792 | 560 | -560 | 0.7071 + 0.7071 i |
| -1352 | -560 | -560 | -1352 | 0.7071 + 0.7071 i |

As expected the ratio (**n**(0) - i**n**(1))/(**n**(2) - i**n**(3)) converges to the largest (absolute value) root, namely $(1 + i)/\sqrt{2}$ of polynomial p(x) in (20). Choosing **R'** = - i**I** + **R** would yield the other root, namely $-(1 + i)/\sqrt{2}$.

## 6. Conclusion

In our earlier papers [2,3,4,5], we had adapted the method of matrix iteration for finding some of the real roots of an integer coefficient polynomial using replacement sequences. The method, though impractical, is interesting because it requires only operations of replacement of sequences and counting of symbols. In this paper, we have shown that the method can be extended to obtain complex roots also. Only one division is needed at the last stage. Even this operation can be carried out geometrically, using a ruler and compass only.